\documentclass[11pt]{elsarticle}
\date{}

\usepackage{graphics} 
\usepackage{subcaption}
\usepackage{graphicx}
\usepackage{epsfig} 
\usepackage{times} 
\usepackage{amsmath} 
\usepackage{mathtools}
\usepackage{amssymb} 
\usepackage{bbold}
\usepackage{color} 
\usepackage{url}
\usepackage{bm}
\usepackage[acronym]{glossaries}

\usepackage{enumitem}
\newlist{lyap}{enumerate}{1}
\setlist[lyap]{label = {\bf (L\arabic*)}, resume}
\newlist{cons}{enumerate}{1}
\setlist[cons]{label = {\bf (C\arabic*)}, resume}

\newcommand{\bnull}{{\boldsymbol 0}}
\newcommand{\bff}{{\bf f}}

\newcommand{\bR}{{\bf R}}

\newcommand{\be}{{\mathbf e}}

\newcommand{\bz}{{\bf z}}
\newcommand{\bx}{{\bf x}}
\newcommand{\bFF}{{\bf F}}

\newcommand{\fS}{\mathfrak S}

\newcommand{\bv}{{\bf v}}

\newcommand{\bphi}{{\boldsymbol \phi}}

\newcommand{\ul}[1]{\underline{#1}}
\newcommand{\ol}[1]{\overline{#1}}

\newcommand{\co}{\operatorname{co}}

\newcommand{\n}{\nonumber}

\newcommand{\cC}{{\mathcal C}}
\newcommand{\cD}{{\mathcal D}}

\newcommand{\cH}{{\mathcal H}}

\newcommand{\cJ}{{\mathcal J}}
\newcommand{\cK}{{\mathcal K}}

\newcommand{\cP}{{\mathcal P}}

\newcommand{\cT}{{\mathcal T}}

\newcommand{\R}{\mathbb{R}}%
\newcommand{\N}{\mathbb{N}}%

\newcommand{\red}[1]{{\color{red} #1}}

\usepackage{multirow}
\usepackage[acronym]{glossaries}

\newcommand{\fela}[1]{}

\newtheorem{theorem}{Theorem}
\newtheorem{corollary}{Corollary}

\newacronym{lp}{LP}{linear programming}
\newacronym{cpa}{CPA}{continuous and piecewise affine}
\newacronym[longplural={linear matrix inequalities}]{lmi}{LMI}{linear matrix inequality}
\newacronym{qclf}{QCLF}{quadratic common Lyapunov function}
\newacronym{clf}{CLF}{common Lyapunov function}
\newacronym{pd}{p.d.}{positive definite}
\newacronym{nd}{n.d.}{negative definite}

\newcommand{\acr}[1]{\acrshort{#1}}
\newcommand{\acrpl}[1]{\acrshortpl{#1}}
\usepackage{lineno}
\makeatletter
\newcommand*\bigcdot{\mathpalette\bigcdot@{.5}}
\newcommand*\bigcdot@[2]{\mathbin{\vcenter{\hbox{\scalebox{#2}{$\m@th#1\bullet$}}}}}
\makeatother

\begin{document}
\begin{frontmatter}

\title{
Linear Programming based Lower Bounds on Average Dwell-Time via Multiple Lyapunov Functions \tnoteref{This work was supported in part by the Icelandic Research Fund
under Grant 228725-051.}
}

\author[inst1]{Sigurdur Hafstein}

\affiliation[inst1]{organization={University of Iceland, Faculty of Physical Sciences},
            addressline={Dunhagi 5},
            city={Reykjavik},
            postcode={107},
            country={Iceland}}
            
\ead{shafstein@hi.is}

\author[inst2]{Aneel Tanwani}

\affiliation[inst2]{organization={CNRS -- LAAS, University of Toulouse},
            city={Toulouse},
            country={France}}

%
\ead{aneel.tanwani@cnrs.fr}


\thispagestyle{empty}
\pagestyle{empty}

\begin{abstract}
With the objective of developing computational methods for stability analysis of switched systems, we consider the problem of finding the minimal lower bounds on average dwell-time that guarantee global asymptotic stability of the origin. Analytical results in the literature quantifying such lower bounds assume existence of multiple Lyapunov functions that satisfy some inequalities. For our purposes, we formulate an optimization problem that searches for the optimal value of the parameters in those inequalities and includes the computation of the associated Lyapunov functions. In its generality, the problem is nonconvex and difficult to solve numerically, so we fix some parameters which results in a linear program (LP). For linear vector fields described by Hurwitz matrices, we prove that such programs are feasible and the resulting solution provides a lower bound on the average dwell-time for exponential stability. Through some experiments, we compare our results with the bounds obtained from other methods in the literature and we report some improvements in the results obtained using our method.
\end{abstract}

\end{frontmatter}

\thispagestyle{empty}

\section{Introduction}

Switched systems comprise a family of dynamical subsystems orchestrated by a switching signal that activates one of these  subsystems at a given time. This abstract framework has been useful in modeling a class of hybrid systems with continuous and discrete dynamics. Another common source of switched systems is uncertainty quantification in continuous-time systems and the associated differential inclusions. Stability analysis of switched systems, therefore, has gathered a lot of attention in the literature. The references \cite{Liberzon2003switched,SWMWK2007revhybrid} provide a comprehensive overview of the different approaches on this topic.

When analyzing stability under arbitrary switching, existence of a common Lyapunov function is a necessary and sufficient condition for the asymptotic stability of an equilibrium of the switched system \cite{DayaMart99}. Thus, over the years, a lot of attention in the literature has been given to computing a common Lyapunov function for the switched system under different hypotheses. For some results in this direction, the reader may refer to \cite{ahmadi14} for  discrete-time systems, and \cite{DellPasq22} for continuous-time systems.
Particularly relevant to this paper is the technique based on the construction of \acrfull{cpa} Lyapunov functions, which is reviewed in~\cite{GiHa2015review}. The papers \cite{Baier2012,Ha2022slidingmodesDI} present the adaptation of computing \acr{cpa} Lyapunov functions in case of arbitrarily switching systems. However, such methods have not yet been used in the context of constrained, or dwell-time based, switched systems.

For certain applications, existence of common Lyapunov function is a stringent requirement, and may not hold for the given system data. For that reason, when the individual subsystems are asymptotically stable and one can not compute a common Lyapunov function, it is natural to ask how we can guarantee stability for a certain class of switching signals (which is smaller than the set of switching signals with arbitrary switching). The works \cite{Morse96} and \cite{HespMorse99}  studied the  stability of switched systems by putting a bound on how fast the switches can occur. Depending on the system data, lower bounds were derived on the (average) dwell-time which ensures global asymptotic stability if the length of interval between two consecutive switches (on average) is greater than the derived lower bound. A tutorial like exposition of these concepts also appears in \cite[Chapter~3]{Liberzon2003switched}. Several works have followed up to extend this idea in several directions.
Some generalizations have been addressed in the recent papers \cite{LiuTanLib22, DellTanw22} with nonlinearities in the system data.

Computational methods with multiple Lyapunov functions for getting best possible lower bounds on the dwell-time have not received much attention in the literature. The references~\cite{Bri15,GerCol06,LMC2021css,Mitra06} provide some algorithms for calculating lower bounds on the dwell-time in the linear case. Among these, the papers \cite{LMC2021css, Mitra06} build on dwell-time bounds obtained from multiple Lyapunov functions, which is also the case for this article. The authors of \cite{Mitra06} developed optimization-based methods for the automatic verification of dwell-time properties. On the other hand, \cite{LMC2021css} proposes some relaxations in the form of sequential convex programs to compute lower bounds on the average dwell-time.
With similar motivation, this article studies computational methods for computing best possible lower bounds on the (average) dwell-time using \acrfull{lp} methods. In fact, our approach uses techniques based on the construction of CPA Lyapunov functions, under the constraints that are normally imposed for dwell-time based stability conditions. For a given family of dynamical subsystems with asymptotically stable origin, the question of interest is to find the {\em smallest} lower bound on the  dwell-time, which ensures asymptotic stability of the switched system under the so-called compatibility constraints. Such questions can be formulated as an optimization problem and in its full generality, it is a nonconvex problem, even when dealing with linear subsystems and quadratic Lyapunov functions for individual subsystems.

In this paper, we provide a new technique for solving the optimization problem that corresponds to the computation of a minimum average dwell-time that ensures stability. The intermediate step in getting this bound is to first compute the Lyapunov functions for individual subsystems satisfying certain inequalities. In our work, we search for these Lyapunov functions from the family of continuous piecewise affine functions, in contrast to the quadratic ones. This is done by discretizing the state space into simplices and solving for the values of the Lyapunov functions at the vertices of the simplices, using some inequality constraints. The resulting optimization problem actually turns out to be a linear program. The solution to this linear program provides us with a Lyapunov function for each subsystem and also a dwell-time bound.

The remainder of the paper is organized as follows: we recall some basic results on (average) dwell-time stability in Section~\ref{sec:setup} and describe the problem being studied in this paper. Section~\ref{sec:lp} provides the \acr{lp} formulation of the proposed problem along with some results about the feasibility of these programs. We provide some simulations and comparisons with other methods in Section~\ref{sec:simuls}, followed by  concluding remarks in Section~\ref{sec:conclusions}.

\section{Problem Setup}\label{sec:setup}

We consider time-dependent switched dynamical systems described as
\begin{equation}\label{eq:swSys}
\dot \bx = \bff_\sigma(\bx)
\end{equation}
where, for some given index set $\cP \subset \N$, the function $\sigma:[0,\infty) \to \cP$ is piecewise constant and right-continuous, called the {\em switching signal}. The discontinuities of $\sigma$, called the {\em switching times}, are assumed to be locally finite. The vector fields $\bff_i:\R^n \to \R^n$, for each $i \in \cP,$ are assumed to be locally Lipschitz and with $\bff_i({\boldsymbol 0}) = {\boldsymbol 0}$. We say that a switching signal $\sigma$ has an average dwell-time $\tau_a > 0$, if there exists $N_0 > 0$, such that
\[
N_\sigma(t,s) \le N_0 + \frac{t-s}{\tau_a},
\]
where $N_\sigma(t,s)$ denotes the number of switches over the interval $(s,t)$. The set of all switching signals with average dwell-time $\tau_a$ is denoted by $\Sigma_{\tau_a}$.

For the stability of the origin for such systems, let us recall the following result, which follows from \cite{HespMorse99}, \cite[Chapter~3]{Liberzon2003switched}, and \cite[Theorem~1]{LiuTanLib22}:

\begin{theorem}\label{thm:mainRecall}
Suppose that
there exist $\cC^1$ Lyapunov functions $V_p:\R^n\to \R_{\geq 0}$, $i\in\cP$, satisfying the following:
\begin{lyap}[leftmargin=*]
\item\label{lyap:bounds} There exist $\underline\alpha,\overline\alpha\in\cK_\infty$ such that
\begin{equation}\label{sandwich_sub_sys}
    \underline\alpha(\|\bx\|)\leq V_i(\bx)\leq \overline\alpha(\|\bx\|),\quad \forall \bx\in \R^n, i\in\cP.
\end{equation}

\item\label{lyap:flow} There exists a Lipschitz function ${\rho} \in \cK$, such that, for every $i \in \cP$,
\begin{equation}\label{continuous_sub_sys_stable}
    \nabla V_i(\bx) \bigcdot \bff_i(\bx)  \leq - {\rho}(V_i(\bx)) \quad \forall \bx\in \R^n,
\end{equation}

\item\label{lyap:jump} There exists $\chi\in\cK_\infty$ such that, for every $i, j \in \cP$, 
\begin{equation}\label{discrete_sub_sys}
    V_j(\bx) \leq \chi (V_i(\bx)) \quad\forall \bx\in \R^n, i \neq j.
\end{equation}
\end{lyap}
Then the origin is globally asymptotically stable for the switched system \eqref{eq:swSys}, uniformly over the set $\Sigma_{\tau_a}$, for $\tau_a$ satisfying
\begin{equation}\label{eq:bndTauNonlin}
\tau_a > \sup_{s> 0} \int_s^{\chi(s)}\frac{1}{\rho(r)}dr.
\end{equation}
\end{theorem}

The lower bound on the average dwell-time given in \eqref{eq:bndTauNonlin}, with nonlinear functions, has also appeared in the context of impulsive systems in \cite{SamoPere95}, \cite{DashMiron13}.  In what follows, we will restrict ourselves to linear subsystems where the functions $\alpha$ and $\rho$ can be taken as linear. Writing \eqref{eq:bndTauNonlin} for such cases allows us to better understand the degrees of freedom at our disposal for minimizing the lower bounds on $\tau_a$.




\subsection{Corollaries and special cases}
Let us present two corollaries to this result depending on the class of functions chosen for $V_i$ and the vector fields $\bff_i$, $i \in \cP$.

\begin{corollary}
Assume that for each $i \in \cP$, there exist symmetric, positive definite matrices $P_i \succ 0$, such that
\begin{subequations}\label{eq:linBMIs}
\begin{align}
A_i^\top P_i + P_i A_i + \alpha P_i \preceq 0, & \quad \text{ for every } i \in \cP, \\
P_j \preceq \mu P_i, & \quad \text{ for every } j \neq i
\end{align}
\end{subequations}
for some $\alpha > 0$ and $\mu \ge 1$. Then the switched system~\eqref{eq:swSys} with $\bff_i(\bx) = A_i \bx $, $i \in \cP$, has a globally exponentially stable equilibrium at the origin, uniformly over the set $\Sigma_{\tau_a}$, where $\tau_a$ satisfies
\begin{equation}\label{eq:dtBndLin}
\tau_a > \frac{\ln(\mu)}{\alpha}.
\end{equation}
\end{corollary}

\medskip

Thus, for linear dynamics and quadratic Lyapunov functions, we can get a lower bound on the average dwell-time by solving matrix inequalities \eqref{eq:linBMIs}. If we take, $\alpha > 0$, $ \mu \ge 1$, and $P_i \succ 0$ as the unknowns in \eqref{eq:linBMIs}, then these inequalities are not linear with respect to the unknowns, and it is difficult to compute a solution. The article \cite{LMC2021css} addresses the problem of minimizing $\frac{\ln(\mu)}{\alpha}$ subject to inequalities \eqref{eq:linBMIs} by proposing convex relaxations.

For the algorithms proposed in this paper, we first need a corollary to Theorem~\ref{thm:mainRecall} with continuous Lyapunov functions, and norm-like bounds on the growth and Dini-derivative of such functions. For a continuous function $V:\R^n \to \R$, we define the Dini-derivative along the solutions of the system $\dot \bx = \bff_i(\bx)$ as
\[
D^+ V(\bx,\bff_i(\bx)):=\limsup_{h\to 0+}\frac{V(\bphi_i(h,\bx))-V(\bx)}{h}
\]
where $t \mapsto \bphi_i(t,\bx)$ is an absolutely continuous function that satisfies $\bphi_i(0,\bx) = \bx$ and $\dot \bphi_i(t,\bx)= \bff_i(\bphi_i(t,\bx))$ almost surely for $t\ge 0$. We use this notion to state the following corollary to Theorem~\ref{thm:mainRecall}, obtained by taking $\ul \alpha$, $\ol \alpha$, $\rho$ to be homogenous functions, and $\chi$ being linear.
\begin{corollary}\label{cor:cpaDT}
Assume that for each $i \in \cP$, there exist continuous functions $V_i:\R^n \to \R$, such that
\begin{subequations}\label{eq:linProgGen}
\begin{align}
    \underline a \|\bx\|^d \le V_i(\bx) \le \overline{a} \|\bx\|^d, & \quad \text{ for every } i \in \cP, \\
    D^+V_i(\bx,\bff_i(\bx)) \le - \alpha \|\bx\|^d, & \quad \text{ for every } i \in \cP, \\
    V_i(\bx) \le \mu V_j(\bx), & \quad \text { for every } j \neq i.
\end{align}
\end{subequations}
for some $d>0$, $\overline a \ge \underline a > 0$, $\alpha > 0$ and $\mu \ge 1$; Then the switched system $\dot \bx = \bff_\sigma(\bx) $ has a globally exponentially stable equilibrium at the origin, uniformly over the set $\Sigma_{\tau_a}$, with $\tau_a$ satisfying,
\begin{equation}\label{eq:dtBndContV}
\tau_a > \frac{\overline{a}\ln(\mu)}{\alpha}.
\end{equation}
\end{corollary}
\medskip

\noindent
In contrast to Theorem~\ref{thm:mainRecall}, the proof of Corollary~\ref{cor:cpaDT} using Dini derivative requires some care but essentially follows similar concepts. In this paper, we will build on the statement of Corollary~\ref{cor:cpaDT} and, in particular, address the following problem:

\paragraph*{Problem statement} With $d=1$ and for fixed values of $\underline a > 0$, $\overline a > 0$, and $\mu \ge 1$, find piecewise linear functions $V_i$, for each $i \in \cP$, that satisfy \eqref{eq:linProgGen} while maximizing $\alpha$.

The reason for fixing the constants $\underline a > 0$, $\overline a > 0$, and $\mu \ge 1$ is that the foregoing problem then transforms into a linear program. We will provide the formulation of this linear program and discuss its feasibility in the next section. For the sake of clarity  in this conference paper, we 
present our ideas for the linear vector fields but similar concepts can be extended to nonlinear systems.

\subsection{Quadratic functions and matrix inequalities}

Before discussing the \acr{lp} problem and \acr{cpa} Lyapunov functions, let us first look at the inequalities \eqref{eq:linProgGen} for the case $d = 2$ more carefully. In this case, we let $V_i(\bx) = \bx^\top P_i \bx$, with symmetric and positive definite $P_i \in \R^{n \times n}$. In particular,
 \eqref{eq:linProgGen} takes the following form, where $I$ is the identity matrix:
\begin{equation}\label{eq:LMI1}
\begin{cases}
    \underline a I \preceq P_i \preceq \overline a I, & \quad i \in \cP, \\
    A_i^\top P_i+P_i A_i \preceq -\alpha I, & \quad i \in \cP, \\
    P_i \preceq  \mu P_j, & \quad i,j \in \cP.
  \end{cases}
\end{equation}
In \eqref{eq:LMI1}, if we fix $\underline a,\overline a >0$ and $\mu \ge 1$, then the inequalities result in LMIs with unknowns $P_i$, $i \in \cP$, which can be solved to maximize $\alpha$.
Practically one can select $\underline{a}$  small, e.g.~$\underline a=10^{-5}$ as we do in our examples, and then a
large enough $\overline{a}>0$ will ensure that  $\tau_a = \overline{a}\ln(\mu)/\alpha$ can be made  minimal for the given $\mu \ge 1$ by maximizing $\alpha$.  Indeed, assume the conditions \eqref{eq:LMI1} are fulfilled for some positive constants $\underline a= \underline a^*,\overline a = \overline a^*, \alpha=\alpha^*$ and $V_i(\bx)=\bx^\top P^*_i\bx$, $P^*_i\in\R^{n\times n}$ symmetric, $i\in \cP$. Then we have the
lower bound $\tau_a^*= \overline{a}^*\ln(\mu)/\alpha^*$  on the average dwell-time.  
Now fix new constants $\underline{a}, \overline{a}>0$ such that  $\overline{a}/\underline{a} \ge \overline{a}^*/\underline{a}^*$ and 
set $\alpha=(\overline{a}/\overline{a}^*)\alpha^* $ and  $P_i=(\overline{a}/\overline{a}^*)P_i^*$.  Then, for $V_i(\bx):=\bx^\top P_i\bx=(\overline{a}/\overline{a}^*)\bx^\top P^*_i\bx$, we have $V_i(\bx)\le \mu V_j(\bx)$, and
$$
\underline{a}\|\bx\|^2 \le \underline{a}^*\frac{\overline{a}}{\overline{a}^*} \|\bx\|^2 \le V_i(\bx) = (\overline{a}/\overline{a}^*) \bx^\top P^*_i\bx \le \overline{a}\|\bx\|^2
$$
for $i,j \in \cP$.  That is, the constraints \eqref{eq:LMI1} are fulfilled with these values of $\underline{a}, \overline{a},\alpha,P_i$ and further, for the lower bound on the average dwell-time we have
$$
\tau_a=\frac{\overline{a}\ln(\mu)}{\alpha} =  \frac{\overline{a}\ln(\mu)}{\frac{\overline{a}\alpha^*}{\overline{a}^*}}=\frac{\overline{a}^*\ln(\mu)}{\alpha^*}=\tau_a^*.
$$
In other words, for a fixed $\mu \ge 1$, if there is a solution to \eqref{eq:LMI1} for some choice of $\ul a^*$, $\ol a^*$, $\alpha^*$ which yields the bound $\tau_a^*$ for the average dwell-time, then by choosing $\ol a/ \ul a$ large enough, one can always find another solution to \eqref{eq:LMI1} which gives at least as good a bound on average dwell-time as $\tau^*$. Thus, given $\mu \ge 1$, maximizing $\alpha>0$ under the constraints \eqref{eq:LMI1} for a fixed $\underline{a},  \overline{a}>0$ delivers as good lower bounds on the average dwell-time $\tau_a$ as minimizing $\tau_a= \overline{a} \ln(\mu)/\alpha$, where both $\overline{a}$ and $\alpha$ are variables, given that $\overline{a}/\underline{a}$ is large enough.


Note that in the setting of the LMI problem \eqref{eq:LMI1} we are searching for quadratic Lyapunov functions for the individual subsystems $\dot\bx=A_i\bx$, $i \in \cP$, which can be conservative. In the next section we consider a similar approach for modeling the conditions \eqref{eq:linProgGen} using piecewise linear Lyapunov functions and an LP formulation to compute them. Due to the foregoing observation, when solving \eqref{eq:linProgGen} using an LP formulation, we will fix $\ul a $ and $\ol a$ with $\ul a/\ol a$ large enough and maximize $\alpha$.

\section{Continuous Piecewise Affine Lyapunov Functions and Linear Programming Formulation}\label{sec:lp}
Our LP approach to compute piecewise linear Lyapunov functions fulfilling the conditions  \eqref{eq:linProgGen}  is based on the so-called CPA method to compute Lyapunov functions, see e.g.~\cite{Mar2002cpa,Baier2012,GiHa2014CPArev,Ha2022slidingmodesDI}.
Its description is somewhat more involved than the \acr{lmi} approach, because it is based on partitioning a neighborhood of the origin into simplices and the underlying idea behind constructing this collection of simplices, called {\em triangulation}, is described in the next subsection.

\subsection{The Triangulation $\cT_K^{\bFF}$}
Roughly speaking, a triangulation is the subdivision of a subset of $\R^n$ into simplices.
A suitable  concrete triangulation  for our aim of parameterizing Lyapunov functions for the individual subsystems is the triangular-fan  of the triangulation in \cite{GiHa2014ACC}, where its efficient implementation is also discussed.
In its definition, we use the functions $\bR^\cJ:\R^n \to \R^n$, defined for every $\cJ\subset\{1,2,\ldots,n\}$ by
$$
\bR^\cJ(\bx) := \sum_{i=1}^n(-1)^{\mathbb{1}_\cJ(i)}x_i\be_i,
$$
where $\be_i$ is the standard $i$th unit vector in $\R^n$, and ${\mathbb 1}_{\cJ}(i) = 1$ if $i \in \cJ$ and $0$ otherwise.  Thus, $\bR^{\cJ}(\bx)$ is the vector $\bx$, except for a minus has been put in front of the coordinate $x_i$ whenever $i\in \cJ$.

We first define the triangulation $\cT^\text{std}$ and use it to construct the intermediate triangulation $\cT_K$, which in turn is used to define our desired triangulation $\cT_K^{\bFF}$. 

 The standard triangulation $\cT^\text{std}$ consists of the simplices
  $$
  \fS_{\bz\cJ\rho}:= \co\left\{\bx_0^{\bz\cJ\rho},\bx_1^{\bz\cJ\rho},\ldots,\bx_n^{\bz\cJ\rho} \right\},
  $$
  where $\co$ denotes the convex hull, and
  \begin{equation}
  \bx_j^{\bz\cJ\rho} := \bR^\cJ\left(\bz+\sum_{i=1}^j\be_{\rho(i)}\right), \label{xjzJsdef}
  \end{equation}
  for all $\bz\in\N^n_{0}=\{0,1,\ldots\}^n$, all $\cJ\subset \{1,2,\ldots,n\}$, all $\rho\in S_n$, and $j=0,1,\ldots,n$. Here, $S_n$ denotes the set of all permutations of $\{1,2,\ldots,n\}$.

  Now fix a $K\in\N_+=\{1,2,\ldots\}$ and define the hypercube $\cH_K:=[-K,K]^n$.  Consider the simplices $\fS_{\bz\cJ\rho}\subset \cH_K$ in
  $\cT^\text{std}$, that intersect the boundary of $\cH_K$.   We are only interested in those intersections that are $(n-1)$-simplices, i.e.~we take every simplex with vertices
  $\bx_j:=\bR^\cJ\left(\bz+\sum_{i=1}^j\be_{\rho(i)}\right)$, $j\in\{0,1,\ldots,n\}$, where exactly one vertex $\bx_{j^*}$ satisfies $ \|\bx_{j^*}\|_\infty<K$ and the other
  $n$ of the $n+1$ vertices satisfy $\|\bx_j\|_\infty=K$, i.e.~for $j\in\{0,1,\ldots,n\}\setminus \{j^*\}$.
   Then we replace the vertex $\bx_{j^*}$ by $\bnull$; it is not difficult to see that $j^*$ is necessarily equal to $0$.
The collection of such vertices triangulates $\cH_K$ and this new triangulation of $\cH_K$ is our desired triangulation $\cT_K$.

\begin{figure}[!h]
\begin{subfigure}[t]{.45\textwidth}
\centering
\includegraphics[width=.8\linewidth]{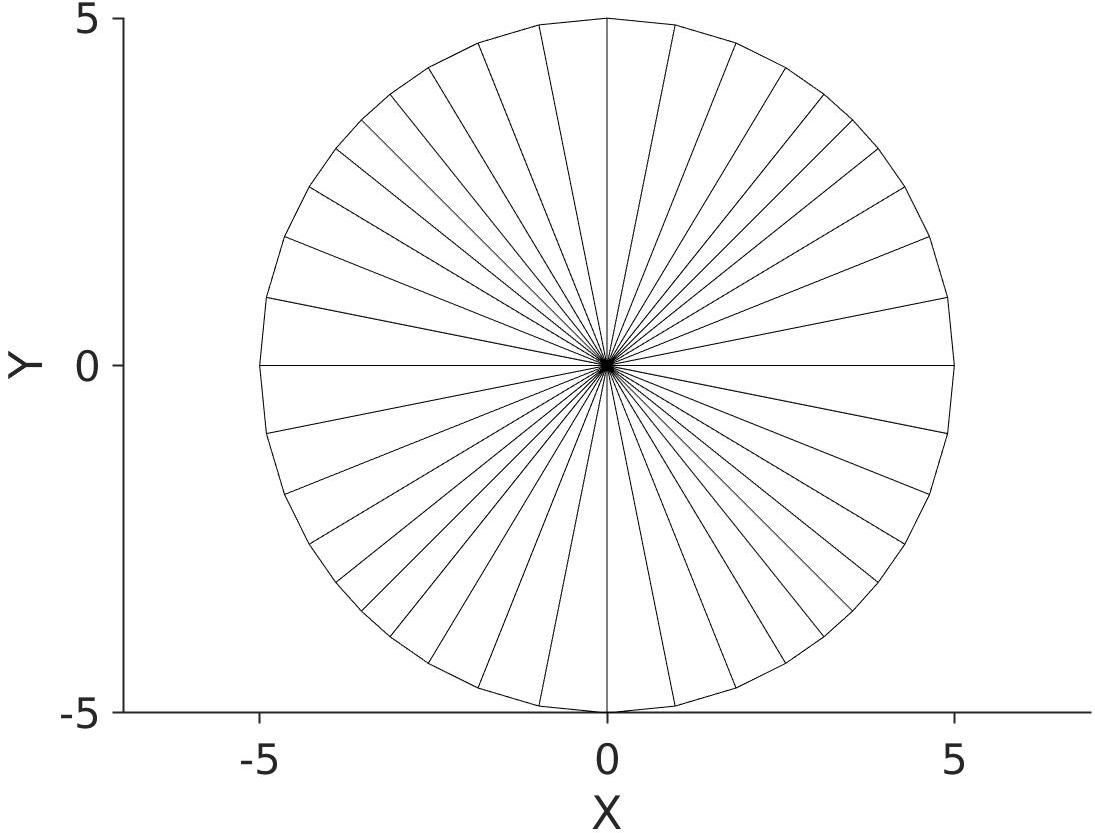}
\caption{The triangulation $\cT_5^{\bFF}$ in two dimensions}
\end{subfigure}
\hfill
\begin{subfigure}[t]{.45\textwidth}
\centering
\includegraphics[width=.8\linewidth]{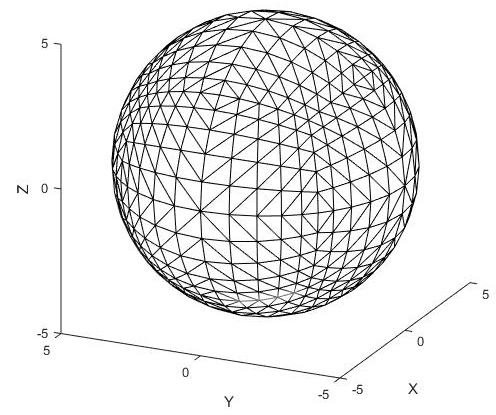}
\caption{The triangulation $\cT_5^{\bFF}$ in three dimension; note that the origin is a vertex of all the tetrahedra in $\cT_5^{\bFF}$.}
\end{subfigure}
\caption{Our proposed triangulation in $\R^2$ and $\R^3$.}
\label{fig:tr3}
\end{figure}


It has been shown \cite{AGGH2019triang} that it is often advantageous in the \acr{cpa} method to map the vertices of the triangulation  by the mapping
$\bFF\colon \R^n\to \R^n$, $\bFF(\bnull)=\bnull$ and
\begin{equation}
\bFF(\bx):= \frac{\|\bx\|}{\|\bx\|_\infty}\bx,\ \ \ \text{if $\bx \neq \bnull$.}
\end{equation}
Note that $\bFF$ maps the hypercubes $\{\bx\in\R^n \colon \|\bx\|_\infty =r\}$ to the spheres $\{\bx\in\R^n \colon \|\bx\| =r\}$.  

Finally, we define the triangulation $\cT_K^{\bFF}$  that will be used in the \acr{lp} problem to parameterize \acr{cpa} Lyapunov functions.
Let $\cT_K^{\bFF}$ be the triangulation consisting of the simplices
$$
  \fS_\nu:= \co\{\bnull,\bFF(\bx_1^\nu),\bFF(\bx^\nu_2),\ldots,\bFF(\bx^\nu_n)\},
  $$
  where 
  $$
\co\left\{\bnull,\bx_1^{\bz\cJ\rho},\bx_2^{\bz\cJ\rho},\ldots,\bx_n^{\bz\cJ\rho} \right\} \in \cT_K.
  $$
The subset of $\mathbb{R}^n$ subdivided into simplices by the triangulation $\cT_K^{\bFF}$ is denoted by
$$
\cD_{\cT_K^{\bFF}}:=\bigcup_{ \fS_\nu \in \cT_K^{\bFF}} \fS_\nu.
$$
Figure \ref{fig:tr3} depicts two exemplary triangulations of the type $\cT_K^{\bFF}$ for two and three dimension with $K=5$.  The implementation of the triangulation is discussed in \cite{Haf2013cpp,peter-siggi-acc}.

\subsection{\acr{lp} Problem
}
\label{sec:LPprob}
We are now ready to state our \acr{lp} problem to parameterize piecewise linear Lyapunov functions for the switched system fulfilling the conditions in \eqref{eq:linProgGen}. For formulating this LP, and showing that its feasibility provides us the lower bound on average dwell-time, we focus our attention on the switched linear systems:
\begin{equation}\label{sys}
\dot \bx = A_\sigma \bx
\end{equation}
with $\sigma:[0,\infty) \to \cP$ being the switching signal, and $A_i \in \R^{n\times n}$, for each $i \in \cP$.


We use three constants $\ul a, \ol a > 0$  and $\mu\ge 1$ in the \acr{lp} problem.  We want the ratio $\ol a/ \ul a$ to be large, as discussed in the last section, and then we want to try out different $\mu\ge 1$ to obtain as good a lower bound on the average dwell-time as possible.

The variables of the \acr{lp} problem are $\alpha\in \R$ and $V_{\bx,i}\in\R$ for every vertex $\bx$ of a simplex in $\cT_K^{\bFF}$ and every $i \in \cP$. 

The objective of the LP problem is to maximize $\alpha$.

The constraints of the \acr{lp} problem are:
\begin{cons}[leftmargin=*]
\item \label{c1}
The first set of constraints is that, for every $i \in \cP$, we set 
$V_{\bnull,i}=0$, and
for every vertex $\bx$ of a simplex in $\cT_K^{\bFF}$ and for every $i \in \cP$: 
\begin{equation}
\label{LP1}
\underline a \|\bx\| \le V_{\bx,i} \le \overline a \|\bx\|
\end{equation}
\item \label{c2}
The second set of constraints is more involved.
For every simplex $\fS_\nu:= \co\{\bnull,\bx_1^\nu,\bx^\nu_2\ldots,\bx^\nu_n\}\in \cT_K^{\bFF}$, we define the matrix $X_\nu=\left(\bx_1^\nu\ \bx_2^\nu\ \cdots \bx_n^\nu\right)$, i.e.~$\bx_k^\nu$ is the $k$th column of $X_\nu$.
Further, we define for every $i \in \cP$, 
the vector of variables $\bv_{\nu,i} =\left(V_{\bx^\nu_1,i}\ V_{\bx^\nu_2,i}\ \cdots \ V_{\bx^\nu_n,i}\right)^\top$.

The constraints are: for every simplex $\fS_\nu\in \cT_K^{\bFF}$, for all  $j=1,\ldots,n$ and all $i \in \cP$: 
\begin{equation}
\label{LP2}
\bv_{\nu,i}^\top X_\nu^{-1} A_i \bx^\nu_j \le -\alpha\|\bx^\nu_j\|.
\end{equation}
Note that these constraints are automatically fulfilled for $j=0$, i.e.~$\bx_j^\nu=\bnull$.
\item\label{c3}
The third set of constraints is:
for every vertex $\bx$ of a simplex in $\cT_K^{\bFF}$ and for every $i,j \in \cP$, 
\begin{equation}
\label{LP3}
V_{\bx,j} \le \mu  V_{\bx,i}.
\end{equation}

\end{cons}

\subsection{Solution to \acr{lp} delivers lower bounds on dwell-time}
In the previous subsection, we formulated an LP which basically specified the constraints in \eqref{eq:linProgGen} at the vertices of the simplices contained in the triangulation. Here we prove that the feasibility of such a program  provides us with piecewise linear Lyapunov functions for the individual subsystems over the entire state space that additionally fulfill (\ref{eq:linProgGen}c), thereby providing a lower bound on the average dwell-time.

Toward this end, assume that the \acr{lp} problem in Section~\ref{sec:LPprob} has a  solution with $\alpha>0$.  We then define
the piecewise linear function $V_i\colon {\cD_{\cT_K^{\bFF}}} \to \R$, for every $i \in \cP$, in the following way:
\begin{itemize}
\item For every $\bx \in \cD_{\cT_K^{\bFF}}$ there exists a simplex $\fS_\nu=\co\{\bnull,\bx_1^\nu,\bx^\nu_2\ldots,\bx^\nu_n\}\in \cT_K^{\bFF}$  such that $\bx\in \fS_\nu$ and there exist a unique
$\lambda\in [0,1]^n$, $\sum_{j=1}^n\lambda_j\le 1$, such that
$\bx=\sum_{j=1}^n\lambda_j \bx^\nu_j$.  We define
\[
V_i(\bx) =\sum_{j=1}^n\lambda_j V_{\bx^\nu_j,i}.
\]
\end{itemize}
It is not difficult to see that the functions $V_i$, $i \in \cP$, are continuous functions that are linear on each simplex $\fS_\nu\in\cT_K^{\bFF}$, in particular each $V_i$ has the constant gradient  $\nabla V_{\nu,i} := \bv_{\nu,i}^\top X_\nu^{-1}$ (row vector) on the interior of $\fS_\nu$, see e.g.~\cite[Rem.~9]{GiHa2014CPArev}.
Hence, for any $\bx \in\fS_\nu\in \cT_K^{\bFF}$, $\bx=\sum_{j=1}^n\lambda_j \bx^\nu_j$, we have for any $i \in \cP$ by \ref{c1} and \ref{c2} 
that
\begin{align}
\nabla V_{\nu,i} \bigcdot A_i \bx&= \bv_{\nu,i}^\top X_\nu^{-1} A_i\sum_{j=1}^n\lambda_j \bx^\nu_j=\sum_{j=1}^n\lambda_j\bv_\nu^\top X_\nu^{-1} A_i\bx^\nu_j \n \\
&\le -\alpha \sum_{j=1}^n\lambda_j \|\bx^\nu_j\| \le -\frac{\alpha}{\ol a} \sum_{j=1}^n\lambda_j V_i(\bx^\nu_j) \n \\
&= -\frac{\alpha}{\ol a} V_i(\sum_{j=1}^n\lambda_j \bx^\nu_j) =  -\frac{\alpha}{\ol a} V(\bx). \label{smalltrick}
\end{align}
Now, for any  $\bx$ in the interior of $\cD_{\cT_K^{\bFF}}$, we have that, for any $i \in \cP$, there exists a simplex $\fS_\nu\in \cT_K^{\bFF}$ and an $h>0$, such that
$$\bx+[0,h]A_i\bx \subset \fS_\nu,
$$
where $\nu$ can depend on both $\bx$ and $i$.
Because $V_i$ is linear on $\fS_\nu$ we have
$$
\limsup_{h\to 0+}\frac{V_i(\bx+hA_{i}\bx)-V_i(\bx)}{h} = \nabla V_{\nu,i} \bigcdot A_i \bx \le -\frac{\alpha}{\ol a} V_i(\bx)
$$
and since this holds true for all $i \in \cP$, we have $D^+V_i(\bx,A_i\bx) \le -\frac{\alpha}{\ol a} V_i(\bx)$.
Since
$$
V_i(\bx)=\sum_{j=1}^n\lambda_j V_{\bx^\nu_j,i} \ge \ul a \sum_{j=1}^n\lambda_j \|\bx^\nu_j\| \ge  \ul a \|\bx\|
$$
 by the constraints \ref{c1}, and for all $i,k \in \cP$,
 $$
 \begin{aligned}
 V_i(\bx)=\sum_{j=0}^n\lambda_j V_{\bx^\nu_j,i} \le \sum_{j=0}^n\lambda_j \mu V_{\bx^\nu_j,k} =\mu V_k(\bx),
\end{aligned}
 $$
  it is clear that the $V_i$ fulfill the constraints \eqref{eq:linProgGen} in the interior of $\cD_{\cT_K^{\bFF}}$.  Just define
 $$
 \ol a':=\max_{\|\bx\|=1}\max_{i=1,2,\ldots,N}V_i(\bx)
 $$
 and we have
 \begin{gather*}
 \ul a \|\bx\| \le V_i(\bx) \le \ol a' \|\bx\|, \\
 D^+V_i(\bx,A_i\bx) \le -\frac{\alpha}{\ol a} V_i(\bx)
 \end{gather*}
 for all $\bx$ in the interior of $\cD_{\cT_K^{\bFF}}$.  Note that we proved \eqref{smalltrick} directly from the constraints and did not go through constraints \eqref{eq:linProgGen} with $\ol a = \ol a'$, which would lead to a worse estimate on $\tau_a$.

 By extending $V_i$ to $\R^n$ in the obvious way, i.e.~for every $\bx\in\R^n$ there exists a $\fS_\nu$ and unique numbers $\lambda_j\ge 0$ such that $\bx=\sum_{j=1}^n\lambda_j \bx^\nu_j$ (a cone defined by the vertices of $\fS_\nu$) and we set $V_i(\bx) =\sum_{j=1}^n\lambda_j V_{\bx^\nu_j,i}$,
 we see that the $V_i$ fulfill the constraints \eqref{eq:linProgGen} on $\R^n$, for each $i \in \cP$.


\section{Simulations}\label{sec:simuls}
We will now test our LP algorithm for several examples from the literature. The class of systems for these simulations is \eqref{sys}. The set $\cP$ and the matrices $\{A_i\}_{i \in \cP}$ will be specified differently for the examples considered here. In the examples, we always fix $\ul a=10^{-5}$ and $\overline a=10$.  We used YALMIP\,/\,sdpt3 and Gurobi to solve the \acr{lmi} and \acr{lp} problems, respectively.

\subsection{Example~1: Dwell-time stable but not under arbitrary switching}
Consider the switched system \eqref{sys}  with
$$
A_1 = \begin{pmatrix}
        -0.1 & -1 \\
        2 & -0.1 \\
      \end{pmatrix}, \ \ \
A_2 = \begin{pmatrix} -0.1 & -2 \\ 1 & -0.1\\ \end{pmatrix}
$$
This example is taken  from \cite[p.~26]{Liberzon2003switched}. It is stable for a certain minimum value for the average dwell-time, but it is not stable under arbitrary switching.
Solving \eqref{eq:LMI1} using the LMI approach, the best $\tau_a$ obtained is $5.1929$ 
with $\mu = 2$.
Using the LP approach, the best $\tau_a$ obtained is $4.5283$ with $\mu = 1.4$ and using $K=500$ in the triangulation.  Using triangulations with fewer triangles delivers a higher lower bound $\tau_a$
for the average dwell time;
$K=50$ gives $\tau_a=5.16493$ with $\mu=1.45$,
$K=100$ gives $\tau_a=4.79315$ with $\mu=1.4$, and
$K=200$ gives $\tau_a=4.62407$ with $\mu=1.4$.  In all cases they are better than the bounds from the LMI approach.
\begin{figure}[!h]
\begin{subfigure}[t]{.24\textwidth}
\centering
\includegraphics[width=.98\linewidth]{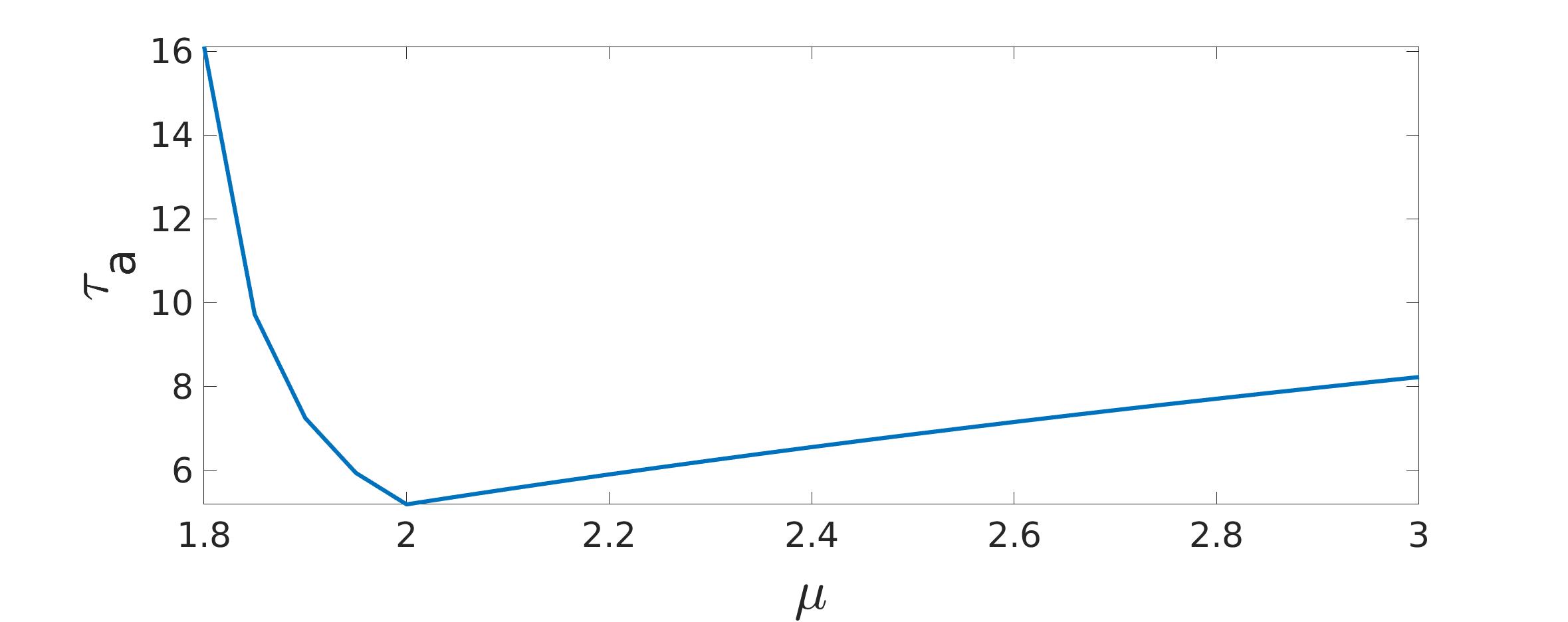}
\caption{LMIs with quadratics}
\end{subfigure}
\begin{subfigure}[t]{.24\textwidth}
\centering
\includegraphics[width=.98\linewidth]{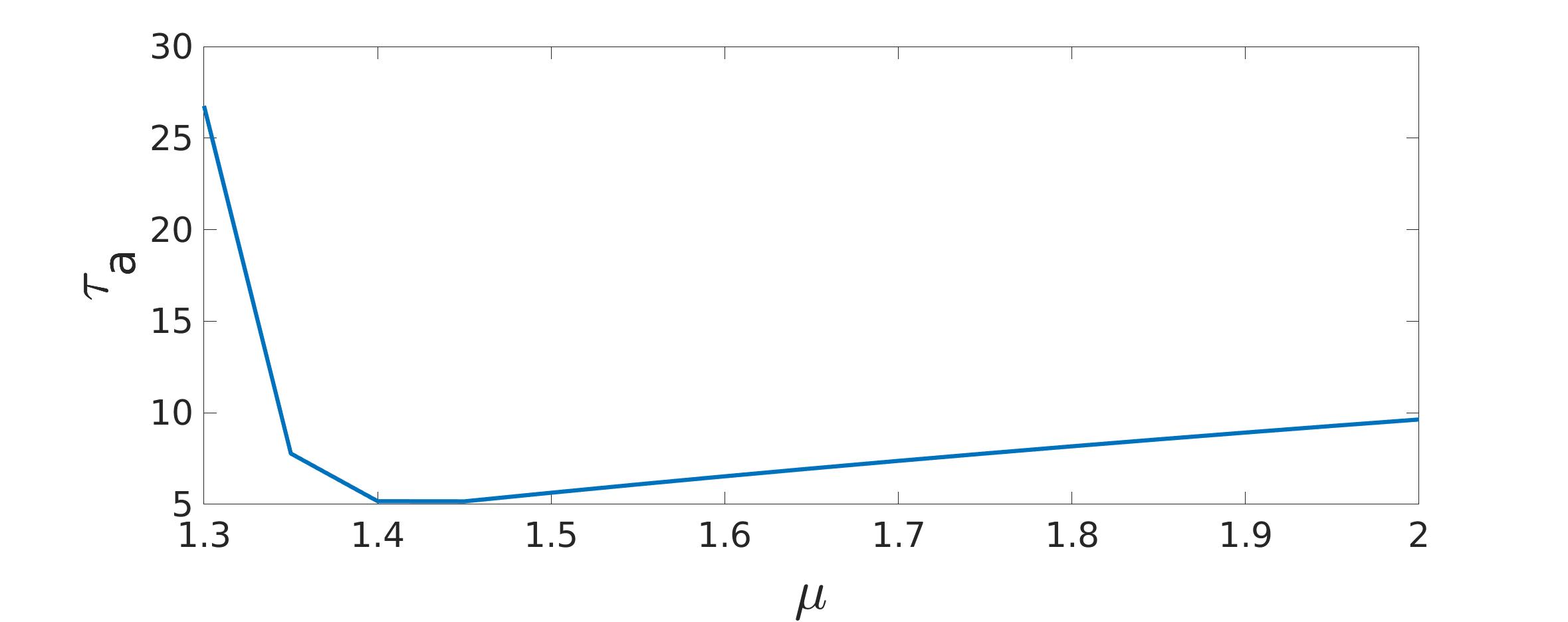}
\caption{LP with $K=50$}
\end{subfigure}
\hfill
\begin{subfigure}[t]{.24\textwidth}
\centering
\includegraphics[width=.98\linewidth]{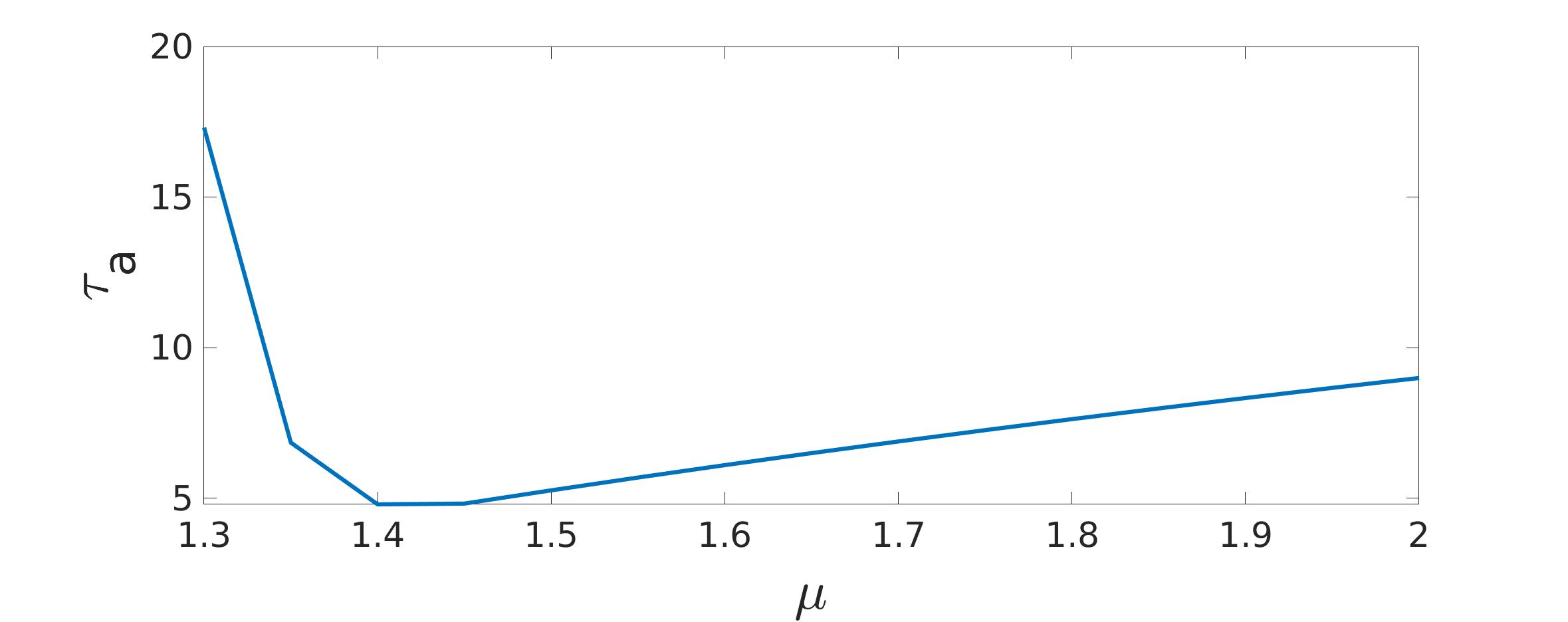}
\caption{LP with $K=100$}
\end{subfigure}
\begin{subfigure}[t]{.24\textwidth}
\centering
\includegraphics[width=.98\linewidth]{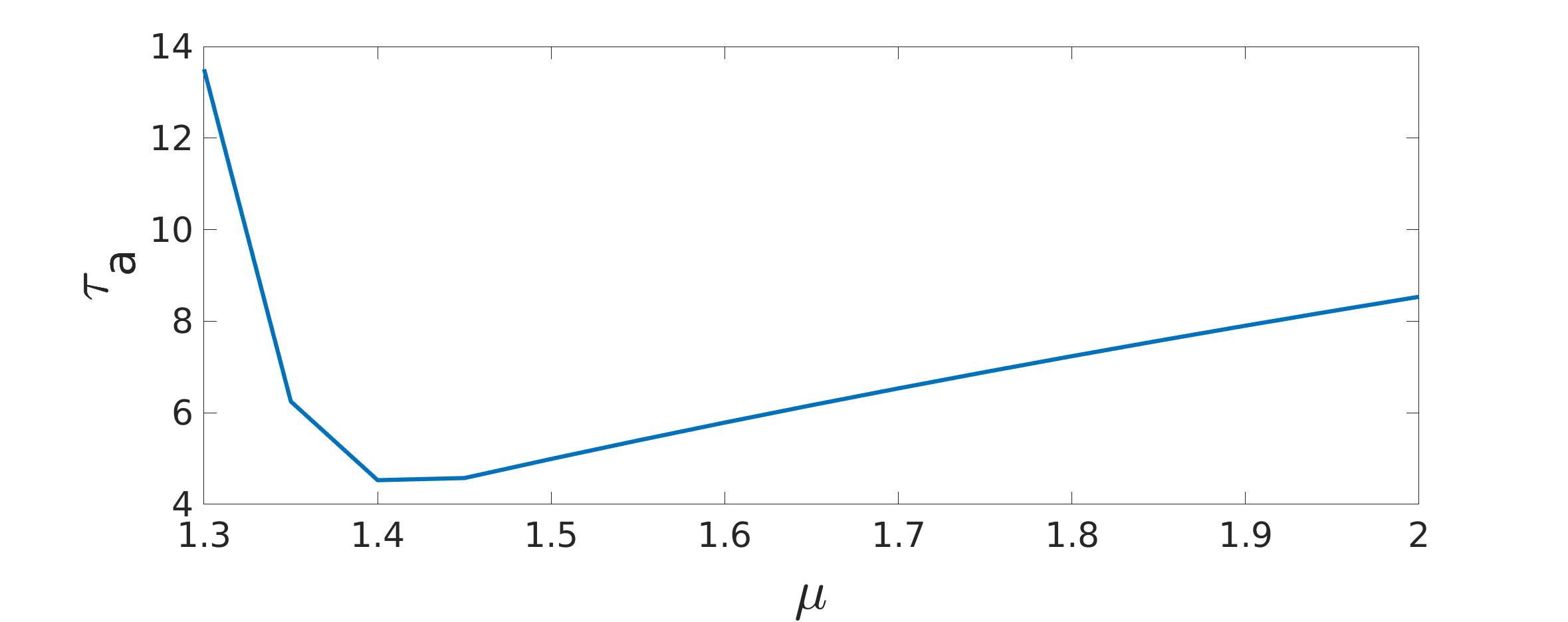}
\caption{LP with $K=500$}
\end{subfigure}
\caption{Plots for $\mu$ versus $\tau_a$ in Example~1.}
\end{figure}

%
%
%
%

\subsection{Example~2: Stable under arbitrary switching but no common quadratic Lyapunov function}
Take $\cP= \{1,2\}$, with the matrices
$$
A_1 = \begin{pmatrix}
        -1 & -1 \\
        1 & -1 \\
      \end{pmatrix}, \ \ \
A_2 = \begin{pmatrix} -1 & -10 \\ 0.1 & -1\\ \end{pmatrix}
$$

This example is taken  from \cite{DayaMart99}. It is stable under arbitrary switching but the matrices $A_1$ and $A_2$ do not share a common quadratic Lyapunov function. This example helps us see the limitation of using LMIs because searching for quadratic certificates in this case is not the best choice.

Solving \eqref{eq:LMI1} using LMIs, the minimum value for $\tau_a$ is $\tau_a= 17.0394$ with $\mu=3.1$. 
Whereas, with our LP approach, $K=20$ gives a solution with $\mu=1$.  Hence, the origin is stable under arbitrary switching ($\tau_a=0$), and we get a common piecewise linear Lyapunov function although no quadratic Lyapunov function exists.

\begin{figure}
\center
\includegraphics[width=0.8\linewidth]{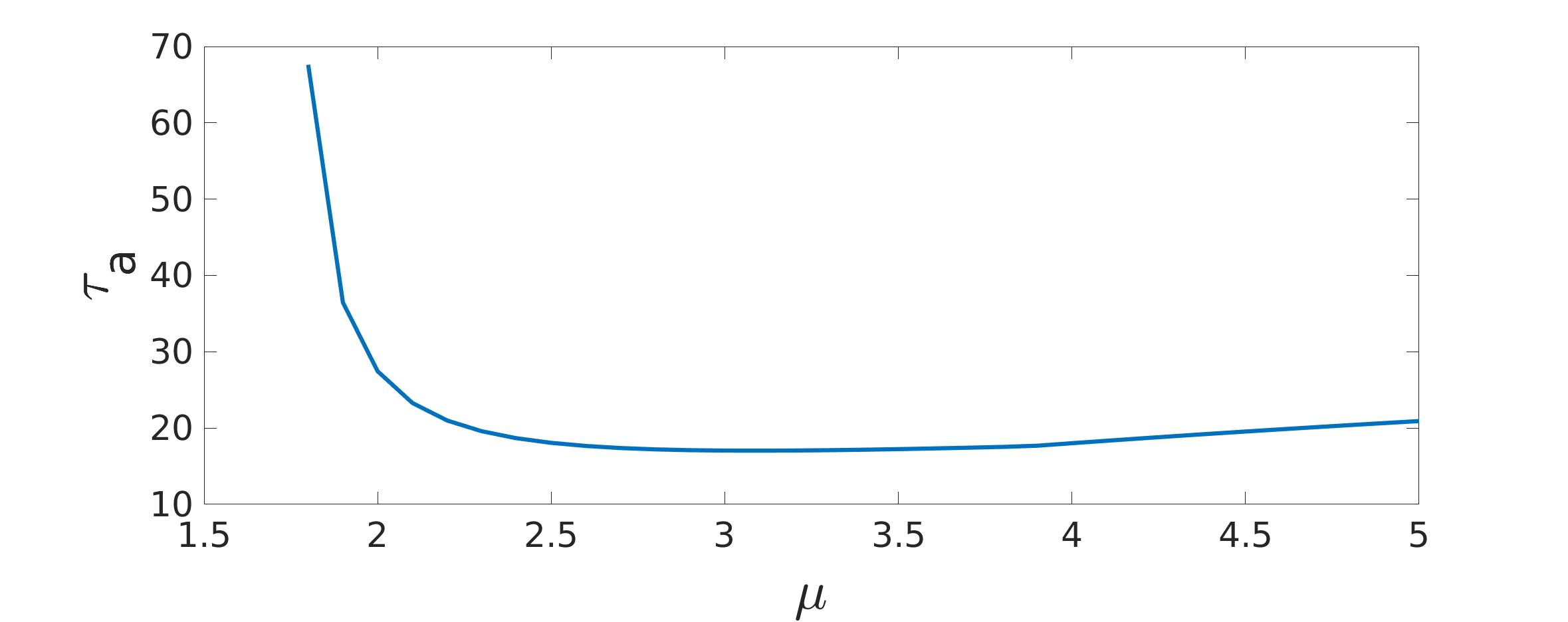}
\caption{LMI based plot for $\mu$ vs $\tau_a$ in Example~2.}\label{fig:X1}
\end{figure}

\subsection{Example 3: Exponentially stable system under arbitrary switching with 5 modes}
We consider the switched system \eqref{sys} with $\cP= \{1,2,3,4,5\}$, where
\begin{align*}
A_1 &= \begin{psmallmatrix} -5&1&2\\0&-5&1\\0&1&-2\\ \end{psmallmatrix},\ \
A_2 = \begin{psmallmatrix}-1&3&1\\0&-2&0\\0&1&-1\\ \end{psmallmatrix}, \\
A_3 &= \begin{psmallmatrix}0&0&3\\-2&-1&-3\\-1&0&-2\\ \end{psmallmatrix},\ \
A_4 = \begin{psmallmatrix}-4&0&-3\\2&-2&4\\1&0&-1\\ \end{psmallmatrix},\\
A_5 &= \begin{psmallmatrix}-1&0&0\\-1&-1&-1\\-3&0&-4\\ \end{psmallmatrix}
\end{align*}
This example was also considered in \cite{LMC2021css} with a graph that determines the switching sequence, and in this particular, we have the star topology.  With $K=6$, we get a solution with $\mu=1$, i.e.~the origin is exponentially stable under arbitrary switching (which is then arbitrary without the graph too).

With the LMI approach in \eqref{eq:LMI1}, the minimum value for the average dwell-time is $\tau_a=4.6870$ 
with $\mu=2.7$.\\


\begin{figure}[!h]
\center
\includegraphics[width=0.8\linewidth]{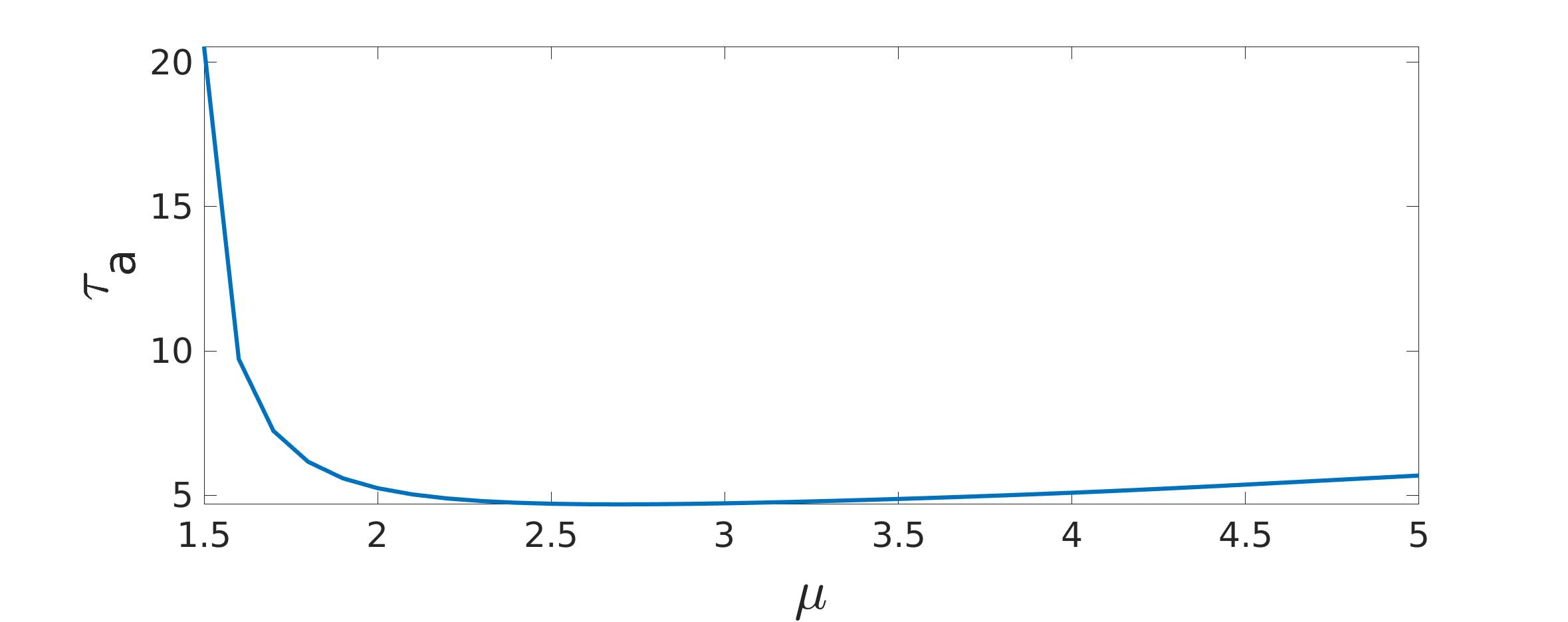}
\caption{LMI based plot for $\mu$ vs $\tau_a$ in Example~3.}\label{fig:X3}
\end{figure}

%

%
%

\section{Conclusions}\label{sec:conclusions}
We considered a \acrfull{lp}  based computational algorithm for computing lower bounds on average dwell-times that ensure asymptotic stability of switched systems. The algorithm is essentially based on gridding the state space into simplices and computing  values for the corresponding  Lyapunov functions at the vertices of these simplices. By choosing appropriate values of the parameters in the inequalities defining the linear program, the solution provides us  lower bounds on the average dwell-time necessary to assure stability. From the simulations, we see in several case studies, that \acr{lp} based bounds are better than the ones based on \acrfullpl{lmi}. This is not really surprising since \acrpl{lmi} restrict the Lyapunov functions to be quadratic, whereas the proposed \acrpl{lp} can potentially approximate a broader class of Lyapunov function templates.
Computing dwell-time via inequalities in \eqref{eq:linProgGen} introduces some conservatism because we optimize over a single parameter $\alpha$ while keeping $\mu$ fixed. The same conservatism is observed in going from \eqref{eq:linBMIs} to \eqref{eq:LMI1}. As a topic of ongoing investigation, we are working out algorithm to optimize $\alpha$ and $\mu$ simultaneously directly using an \acr{lp} version of \eqref{eq:linBMIs}.  Other than understanding the complexity of the  proposed algorithm, we also aim to study the extensions of the proposed algorithm in different directions, which includes the study of generalized lower bounds on average dwell-time \cite{KunDebLib16,LiuTanLib22}.

\bibliographystyle{plain}
\bibliography{LyapBibUniform,biblio}

\end{document}